\pdfoutput=1
\documentclass[a4paper]{article}
\usepackage{a4wide}
\usepackage{microtype}

\usepackage[
  algo2e,
  lined,
  boxed,
  commentsnumbered,
  linesnumbered,
  ruled
]{algorithm2e}

\usepackage{graphicx}
\usepackage{amsmath}
\usepackage{amsfonts}
\usepackage{amssymb}
\usepackage{amsthm}
\usepackage{booktabs}
\usepackage{placeins}
\usepackage{pgf}
\usepackage{tikz}
\usepackage{pgfplots}
\usepackage{subcaption}
\usetikzlibrary{plotmarks}

\newcommand{\R}{\mathbb{R}}

\newcommand{\funspace}{V}
\newcommand{\redspace}{{{}\widetilde{\funspace}}}
\newcommand{\etaspace}{{\funspace_\eta}}
\newcommand{\hdim}{{N}}
\newcommand{\solution}{u_\mu}
\newcommand{\rbsolution}{{{}\widetilde{u}_{\mu}}}
\newcommand{\rbdim}{{\widetilde{\hdim}}}
\newcommand{\etadim}{{\hdim_\eta}}
\newcommand{\etabase}{\psi^\eta}
\newcommand{\Parspace}{\mathcal{P}}
\newcommand{\Residual}{\mathcal{R}_\mu}
\newcommand{\Riesz}{\mathsf{R}}
\newcommand{\norm}[1]{{\left\lVert{#1}\right\rVert}}
\newcommand{\coef}[1]{{\mathsf{#1}}}

\DeclareMathOperator{\spn}{span}

\newtheorem{theorem}{Theorem}[section]

\begin{document}
\title{A numerically stable a posteriori error estimator 
for reduced basis approximations of elliptic equations}
\date{March 17, 2014}

\author{Andreas Buhr%
\thanks{Institute for Computational und Applied Mathematics
Einsteinstrasse 62, 48149 M\"unster, Germany,
\{andreas.buhr, christian.engwer, mario.ohlberger, stephan.rave\}@uni-muenster.de
}
\and Christian Engwer\footnotemark[1]
\and Mario Ohlberger\footnotemark[1]
\and Stephan Rave\footnotemark[1]
}

\maketitle

\begin{abstract}
The Reduced Basis (RB) method is a well established method
for the model order reduction of problems formulated as parametrized
partial differential equations.
One crucial requirement for the application of RB schemes is
the availability of an a posteriori error estimator to reliably estimate the error introduced by the reduction process.
However, straightforward implementations
of standard residual based estimators show poor numerical stability, 
rendering them unusable
if high accuracy is required.
In this work we propose a new algorithm based on representing
the residual with respect to a dedicated orthonormal basis, which is both
easy to implement and requires little
additional computational overhead. A numerical example
is given to demonstrate the performance of the proposed algorithm.
\end{abstract}

\section{INTRODUCTION}
Many problems in science and engineering require the solution
of partial differential equations on large computational domains
or very fine meshes. Even on modern hardware, 
standard discretization techniques for solving these problems can
require many hours or even days of computation, which makes these
approaches inapplicable for many-query situations like, e.g.,
design optimization, where the same problem has to be solved many
times for different sets of parameters.

The Reduced Basis Method (RB) is by now a well-established tool for
the model order reduction of problems formulated as parametrized
partial differential equations.
For a general introduction we refer to \cite{RBBook} and \cite {Haasdonk}.
In an ``offline phase'', a given
high-dimensional discretization is solved for appropriately selected
parameters and a reduced subspace is constructed as the span
of these solution snapshots. In a later ``online phase'', the problem
can be solved efficiently for arbitrary new parameters via Galerkin
projection onto the precomputed reduced space.

One crucial ingredient for the application of RB schemes is
the availability of a quickly evaluable a posteriori error estimator to reliably estimate
the error introduced by the reduction process. Such an estimator is also required
by the weak greedy algorithm, which has been shown to be
optimal for the generation of the reduced spaces \cite{Binev}, to
efficiently perform an exhaustive search of the parameter space for
parameters maximising the reduction error.

For affinely decomposed elliptic problems, a residual based error
estimator is widely used \cite[sec.\ 4.3]{RBBook}.
In order to ensure quick evaluation of the dual norm of the residual, the computation
is decomposed into high-dimensional operations during the ``offline phase'' 
and fast low-dimensional computations during the ``online phase''.
However, as observed by several authors
\cite[pp.\ 148--149]{RBBook}\cite{YANO}\cite{BennerHess}, 
the implementation of this offline/online splitting shows poor numerical
accuracy due to round-off errors which can render the estimator unusable
when the given problem is badly conditioned and high accuracy
is required.
Observations suggest that the estimator typically 
stagnates at a relative error of order $\sqrt{\varepsilon}$,
where $\varepsilon$ is the machine accuracy of the
floating point hardware used.

In the following, we propose a new algorithm to evaluate
the norm of the residual which does not
suffer the severe numerical problems of the traditional approach,
is free of approximations, has only small computational overhead
and is easy to implement.

To our knowledge, 
there is only one other contribution in which a numerically stable algorithm
for evaluation of the estimator is presented
\cite{Casenave12,Casenave13}. 
This approach however comes at the price of
a computationally more expensive ``online phase'' (in \cite{Casenave12})
or increased complexity of
offline computations (in \cite{Casenave13}) by application
of the empirical interpolation method, which in turn requires additional stabilization.
Moreover, a proof for the reliability of the modified estimator is
missing in~\cite{Casenave13}.

The remainder of this
paper is organized as follows: In Section \ref{sec:problem} we introduce
the high-dimensional discrete problem that we will consider in this work. In Section \ref{sec:rbmethod}
we summarize the Reduced Basis method including
the weak greedy algorithm for basis generation.
In Section \ref{sec:errorestimator} we present the
residual based error estimator under consideration, the traditional
algorithm for its evaluation as well as our proposed new algorithm.
Finally, in Section \ref{sec:numericalexample} we 
give a numerical example underlining the improved stability
of our new algorithm.

\section{HIGH-DIMENSIONAL PROBLEM}
\label{sec:problem}
We consider a discrete parametrized elliptic problem of the following form: let
$\funspace$ be a Hilbert space of finite dimension $N$, $f_\mu \in \funspace'$ 
a parametrized linear functional and $a_\mu : \funspace \times \funspace \rightarrow \R$
a parametrized bilinear form such that for an $\alpha_\mu > 0$ we have
\begin{equation}
\label{eq:coercivity}
\alpha_\mu \norm{\varphi}_\funspace^2 \leq a_\mu(\varphi,\varphi) \qquad \forall
\varphi \in \funspace.
\end{equation}
We then search for the solution $\solution \in \funspace$ satisfying
\begin{equation}
\label{eq:mainequation}
a_\mu(\solution,\varphi) = f_\mu(\varphi) \qquad \forall \varphi \in \funspace.
\end{equation}
Note that the existence of a solution follows from the coercivity (\ref{eq:coercivity}) of $a_\mu$ and the finite dimensionality of $\funspace$.
The parameter $\mu$ is confined to be an element of a fixed compact parameter space $\Parspace \subset \R^P$.
Moreover, we assume that $a_\mu$ and $f_\mu$ exhibit
an affine parameter dependence, i.e.~there exist parameter independent
bilinear forms $a^q : \funspace \times \funspace \rightarrow \R$ $(1 \leq q \leq
Q_a)$, linear functionals $f^q \in \funspace'$ $(1 \leq q \leq Q_f)$
 and
coefficient functionals $\theta_a^q : \Parspace \rightarrow \R$ and $\theta_f^q : \Parspace \rightarrow \R$ 
such that
\begin{equation}
\label{eq:affineaf}
a_\mu(\varphi_1,\varphi_2) = \sum_{q=1}^{Q_a} \theta_a^q(\mu) a^q(\varphi_1, \varphi_2)
\qquad\mathrm{and} \qquad
f_\mu(\varphi) = \sum_{q=1}^{Q_f} \theta_f^q(\mu) f^q(\varphi).
\end{equation}
\section{REDUCED BASIS APPROXIMATION}
\label{sec:rbmethod}
Given smooth dependence of the solution $\solution$ on the parameter $\mu$, the
dimension of the manifold of all solutions $\{u_\mu\,\vert\, \mu \in \Parspace\}$
is bounded by $\dim (\Parspace)$ and, thus, is in general of much lower
dimension than $\funspace$.
The Reduced Basis method exploits this fact
by constructing a low-dimensional linear subspace $\redspace \subset \funspace$
of dimension $\rbdim$
in which the solution manifold can be approximated up to a small error. A
reduced solution $\rbsolution \in \redspace$ is then determined by Galerkin
projection of (\ref{eq:mainequation}) onto $\redspace$, i.e.~by solving
\begin{equation}
\label{eq:mainrbequation}
a_\mu(\rbsolution,\widetilde \varphi) = f_\mu(\widetilde \varphi) \qquad \forall \widetilde \varphi \in \redspace.
\end{equation}
The solvability of (\ref{eq:mainrbequation}) again follows from
(\ref{eq:coercivity}).

The reduced space $\redspace$ is constructed from the linear span 
of solutions
to (\ref{eq:mainequation}) for parameters selected by the following greedy search procedure:
Starting with $\redspace_0 := \{ 0 \} \subset \funspace$,
in each iteration step the reduced problem (\ref{eq:mainrbequation})
is solved and an error estimator is evaluated at all parameters
$\mu$ of a given training set $\mathcal{S}_{train} \subset \Parspace$.
If the maximum estimated error is below a prescribed tolerance $tol$,
the algorithm stops. Otherwise, the high-dimensional problem (\ref{eq:mainequation})
is solved for the parameter $\mu^*_n$ maximising the estimated error
and the reduced space is extended by the obtained solution snapshot: 
$\redspace_{n+1} := \redspace_{n} \oplus \spn\{ u_{\mu^*_n}\}$.
\section{RESIDUAL BASED A POSTERIORI ERROR ESTIMATOR}
\label{sec:errorestimator}
An a posteriori error estimator provides a computable upper bound for the
model reduction error $\norm{\solution - \rbsolution}_\funspace$. 
We consider here a widely used error estimator based on the discrete residual
$\Residual \in \funspace^\prime$ given by $\Residual(\rbsolution)(\varphi) := f_\mu(\varphi)
- a_\mu(\rbsolution,\varphi)$ for $\varphi \in \funspace$.

\begin{theorem}
\emph{(Error bound)} 
The model reduction error $\norm{\solution - \rbsolution}_\funspace$
can be bounded using the dual norm of the residual and the coercivity constant of the
bilinear form:
\begin{equation}
\norm{\solution - \rbsolution}_\funspace \leq \frac{1}{\alpha_\mu} \norm{\Residual(\rbsolution)}_{\funspace'}
\end{equation}
\end{theorem}
\begin{proof}
See \cite[eq. 4.28]{RBBook}.
\end{proof}
To calculate the dual norm of the residual $\Residual(\rbsolution)$ we make use of the
fact that the norm of an element of $\redspace^\prime$ is equal to the norm
of its Riesz representative. Denoting by $\Riesz : \funspace' \rightarrow \funspace$
the Riesz isomorphism and assuming the existence of a computable lower bound $\alpha_{\mu,LB} \leq \alpha_\mu$ for the coercivity constant,
we obtain a bound for the error containing only computable quantities:
\begin{equation}
\label{eq:computable_bound}
\norm{\solution - \rbsolution}_\funspace \leq \frac{1}{\alpha_{\mu,LB}}
\norm{\Riesz(\Residual(\rbsolution))}_{\funspace}
\end{equation}
Direct evaluation
of this error bound comprises the calculation of the Riesz
representative and the computation of its norm,
which are both high-dimensional operations. 
However, for the application of the RB method in many-query and real-time situations, 
it is crucial 
that the time for evaluating the a posteriori error estimator in the online phase
is independent of the dimension of $\funspace$. This is also required to make
the use of large parameter training sets $\mathcal{S}_{train}$ feasible, which
is necessary to ensure optimal selection of the snapshot parameter $\mu^*$.
\subsection{Traditional offline/online splitting}

In order to avoid high-dimensional calculations
during the online phase, the residual $\Residual(\rbsolution)$ can be
rewritten using the affine decompositions (\ref{eq:affineaf}) and
a basis representation of $\rbsolution$. Let
$\{\widetilde{\psi}_1, \dots, \widetilde{\psi}_\rbdim\}$
be a basis of $\redspace$ and let  
$
\rbsolution = \sum_{i=1}^{\rbdim} \coef{\rbsolution}_i \widetilde{\psi}_i
$,
then the Riesz representative of the residual is given as
\begin{equation}
	\Riesz(\Residual(\rbsolution)) = 
\sum_{q=1}^{Q_f} \theta_f^q(\mu) \Riesz(f^q) - \sum_{q=1}^{Q_a} \sum_{i=1}^{\rbdim}
\theta_a^q(\mu) \coef{\rbsolution}_i \Riesz(a^q(\widetilde{\psi}_i,\, \cdot\,))\,.
\label{eq:residual}
\end{equation}
To simplify notation, we rename the $\etadim := Q_f$ + $Q_a  \rbdim$ linear 
coefficients $\theta_f^q(\mu)$ and $\theta_a^q(\mu) \coef{\rbsolution}_i$ to $\alpha_k$ 
and the vectors $\Riesz(f^q)$ and $\Riesz(a^q(\widetilde{\psi}_l,\cdot))$ to
$\eta_k$, i.e.~$\Riesz(\Residual(\rbsolution)) = \sum_{k=1}^{\etadim} \alpha_k \eta_k$.
The space $\spn\{\eta_1, \dots, \eta_\etadim\}$ is denoted by $\etaspace$.
For the norm of the residual we obtain
\begin{equation}
\label{eq:traditional_offline_online}
\norm{\Riesz(\Residual(\rbsolution))}_\funspace = \left( \sum_{k=1}^{\etadim}
	\sum_{l=1}^{\etadim} \alpha_k \alpha_l \left( \eta_k , \eta_l
	\right)_\funspace \right)^\frac{1}{2}.
\end{equation}
Using this representation, an offline/online decomposition of the
error bound is possible by pre-computing the inner products $(\eta_k, \eta_l)_\funspace$
during the offline stage. In the
online stage, only the sum in (\ref{eq:traditional_offline_online}) has to be evaluated. As the number of summands is independent of the dimension of $\funspace$,
an online run-time independent of the dimension of $\funspace$ is achieved.

While this approach leads to an efficient computation of the residual
norm, it shows poor numerical stability: in the sum
(\ref{eq:traditional_offline_online}), terms with a relative error of order of machine
accuracy $\varepsilon$ are added. Therefore, the sum shows an absolute error of at least $\varepsilon$ times the largest value of $|\alpha_k \alpha_l (\eta_k, \eta_l)_\funspace|$,
and the error in the norm of the residual 
is thus at least of order $\sqrt{\varepsilon}\cdot\sqrt{\max_{k,l}(|\alpha_k \alpha_l (\eta_k,
	\eta_l)_\funspace|)}$. This is in agreement with the observation that this
algorithm stops converging at relative errors of order
$\sqrt{\varepsilon}$ (see Section \ref{sec:numericalexample}).
\subsection{Improved offline/online splitting}
While the floating point evaluation of (\ref{eq:traditional_offline_online}) shows poor
numerical accuracy, note that the evaluation of 
\begin{equation}
\label{eq:hdevaluation}
\norm{\Riesz(\Residual(\rbsolution))}_\funspace = 
\left(\sum_{k=1}^{\etadim} \alpha_k \eta_k,\sum_{k=1}^{\etadim} \alpha_k \eta_k\right)_{\funspace}^\frac{1}{2}
\end{equation}
is numerically stable.
Based on this observation, we propose a new algorithm to evaluate
$\norm{\Riesz(\Residual(\rbsolution))}_\funspace$ which is offline/online decomposable while maintaining the algorithmic
structure of (\ref{eq:hdevaluation}) to ensure stability.

The algorithm we propose evaluates (\ref{eq:hdevaluation}) in the subspace
$\etaspace$ 
using an orthonormal basis for this space. It comprises three steps:
1. The construction of an orthonormal basis $\Psi^\eta = \{ \etabase_1, \dots,
\etabase_\etadim \}$ of $\etaspace$,
2. the evaluation of the basis coefficients of $\eta_k$ w.r.t.~the basis $\Psi^\eta$ and
3. the evaluation of (\ref{eq:hdevaluation}) using this basis representation.
Note that this approach is offline/online decomposable: Steps 1 and 2 can be 
done offline, without knowing the parameter, while step 3 can be performed online. 
The size of the basis $\Psi^\eta$ does not depend on the dimension of $\funspace$.

\begin{algorithm2e}[t]
\KwIn{vectors $v_i$, $i \in 1, \dots, N$}
\KwOut{orthonormal vectors $v_i$}
\For{$i \leftarrow 1, \dots, N$}{
  $v_i \leftarrow v_i / \norm{v_i}_\funspace $\;
  \Repeat{$\mathrm{newnorm}  > 0.1$}{
    \For{$j \leftarrow 1, \dots, (i-1)$}{
      $v_i \leftarrow v_i - (v_i,v_j)_\funspace ~ v_j$\;
    }
    newnorm $\leftarrow \norm{v_i}_\funspace$\;
    $v_i \leftarrow v_i / \mathrm{newnorm} $\;
  }
}
\caption{Gram-Schmidt with re-iteration}
\label{alg:gram-schmidt}
\end{algorithm2e}

In principle, any orthonormalization algorithm 
applied to $\{\eta_k\}_{k=1}^{\etadim}$ can be used for the computation of
the basis $\Psi^\eta$. Note, however, that the algorithm has to compute the basis
with very high numerical accuracy. As an example, the standard modified
Gram-Schmidt algorithm usually fails to deliver the required accuracy. For the
numerical example in Section \ref{sec:numericalexample}, we have chosen an
improved variant of the modified Gram-Schmidt algorithm, where vectors are
re-orthonormalized until a sufficient accuracy is achieved (Algorithm~\ref{alg:gram-schmidt}).

After the basis $\Psi^\eta$ has been constructed using an appropriate orthonormalization algorithm, 
we can compute for each $\eta_k$ ($1 \leq k \leq \etadim)$ basis representations
$\eta_k = \sum_{i=1}^{\etadim} \overline{\eta}_{k,i} \etabase_i$, where
$
\overline{\eta}_{k,i} = \left( \eta_k , \etabase_i \right) _ \funspace
$ due to the orthonormality of $\Psi^\eta$.
The right-hand side of (\ref{eq:hdevaluation}) can then be evaluated as:
\begin{equation}
\label{eq:newestimator}
\norm{\Riesz(\Residual(\rbsolution))}_\funspace
= \left( \sum_{i=1}^\etadim \left( \sum_{k=1}^\etadim \alpha_k \overline \eta_{k,i} \right)^2 \right)^\frac{1}{2},
\end{equation}
which executes in time independent of the dimension of $\funspace$ and is observed to be numerically stable.
\subsection{Run-time complexities}
During the offline phase, both the traditional and the new algorithm have to calculate
all Riesz representatives appearing in (\ref{eq:residual}). This requires the application of the inverse 
of the inner product matrix for $V$, which 
can be computed in complexity $\mathcal{O}(\hdim \log(\hdim))$ 
with appropriate preconditioners.
As there are $\etadim$ Riesz representatives to be calculated, the overall run-time
of this step is of order $\mathcal{O}(\etadim \hdim \log(\hdim))$.
The traditional algorithm proceeds with calculating all
inner products $(\eta_k,\eta_l)_V$ in (\ref{eq:traditional_offline_online}), having a complexity of $\mathcal{O}(\etadim^2 \hdim)$.
Thus the overall complexity of the offline phase for the traditional algorithm is
$\mathcal{O}(\etadim^2 \hdim + \etadim \hdim \log(\hdim))$.

After computing the Riesz representatives in (\ref{eq:residual}),
the improved algorithm generates the orthonormal basis $\Psi^\eta$.
In practice it was observed that at most four re-iterations per vector are required during
orthonormalization with Algorithm \ref{alg:gram-schmidt}. Thus, 
choosing this algorithm for the generation of $\Psi^\eta$ leads
to a run-time complexity of $\mathcal{O}(\etadim^2 \hdim)$ for this step.
The calculation of the $\etadim^2$ basis coefficients
$\overline \eta_{k,i} = \left( \eta_k , \etabase_i \right) _ \funspace$
has again complexity $\mathcal{O}(\etadim^2 \hdim)$,
resulting in a total complexity of the offline phase for the new algorithm of
$\mathcal{O}(\etadim^2 \hdim + \etadim \hdim \log(\hdim))$, as for the traditional algorithm.

During the online phase, the right-hand sides of (\ref{eq:traditional_offline_online}), resp.~(\ref{eq:newestimator}),
are evaluated using the pre-computed quantities $(\eta_k,\eta_l)_V$, resp.~$\overline \eta_{k,i}$.
In both cases, a run-time of $ \mathcal{O}(\etadim^2)$ is required.

\begin{table}[b]
\caption{Run-time complexities of traditional and new algorithm for evalution of the error estimator.}
\label{tab:complexities}
\begin{center}
\setlength{\tabcolsep}{4.5mm}
\begin{tabular}{rcc}
\toprule
stage & offline & online \\
\midrule
traditional  & $\mathcal{O}(\etadim^2 \hdim) + \mathcal{O}(\etadim \hdim \log(\hdim))$ & $\mathcal{O}(\etadim^2)$ \\
new  & $\mathcal{O}(\etadim^2 \hdim) + \mathcal{O}(\etadim^2 \hdim) + \mathcal{O}(\etadim \hdim \log(\hdim))$ & $\mathcal{O}(\etadim^2)$ \\
\bottomrule
\end{tabular}
\end{center}
\end{table}

All in all, both algorithms for evaluating (\ref{eq:computable_bound}) show the same
run-time complexity, in the online phase as well as during the offline phase (Table \ref{tab:complexities}).
Note that 
$\mathcal{O}(\etadim^2) = \mathcal{O}(Q_f^2 + Q_a^2 \rbdim^2) = \mathcal{O}(\rbdim^2)$
for increasing reduced space dimensions.
\section{NUMERICAL RESULTS}
\label{sec:numericalexample}
\begin{figure}
\begin{center}
\includegraphics[scale=0.5, trim=00 230 0 200, clip]{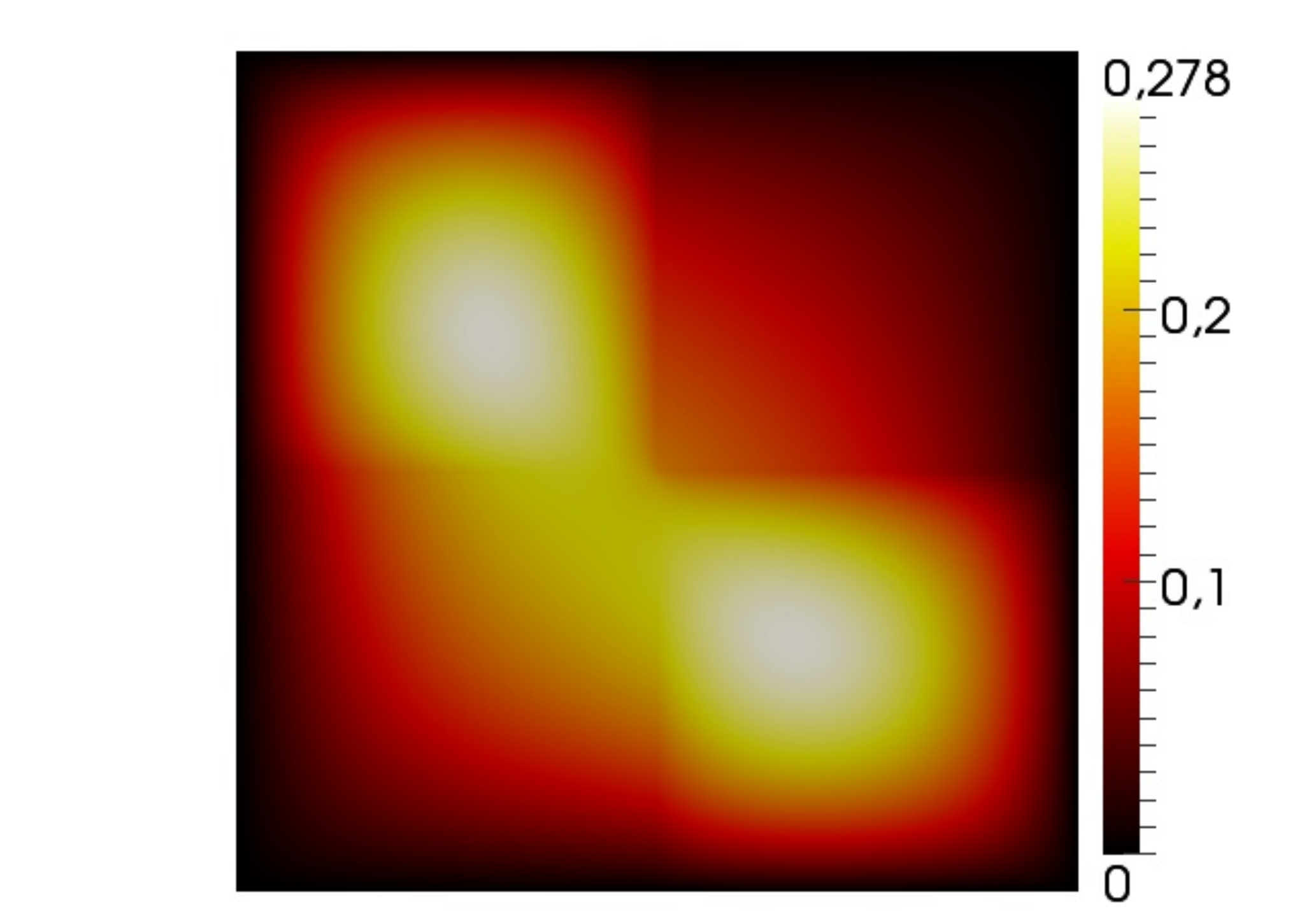}
\end{center}
	\caption{High-dimensional solution of (\ref{eq:thermalblock}) for $\mu = (0.1, 1.0, 0.4,
		1.0)$}
\label{fig:solution}
\end{figure}
In order to verify the improved numerical stability of our proposed algorithm, we
considered an elliptic ``thermal block'' problem on the domain $\Omega = [0, 1]^2$ of the form
\begin{equation}
	\label{eq:thermalblock}
	- \nabla \cdot (\sigma_\mu \nabla u_\mu) = 1, \qquad u_\mu \in H^1_0(\Omega),
\end{equation}
with heat conductivity 
$\sigma_\mu = \sum_{i,j=0}^1 \mu_{ij}\cdot\chi_{[i/2, (i+1)/2]\times [j/2,
	(j+1)/2]}$, 
denoting by $\chi_A$ the characteristic function of the set $A$. The parameters
$\mu=\{\mu_{ij}\}_{i,j=0}^1$ were allowed to vary in the space $\mathcal{P} = [0.1, 1.0]^4$.

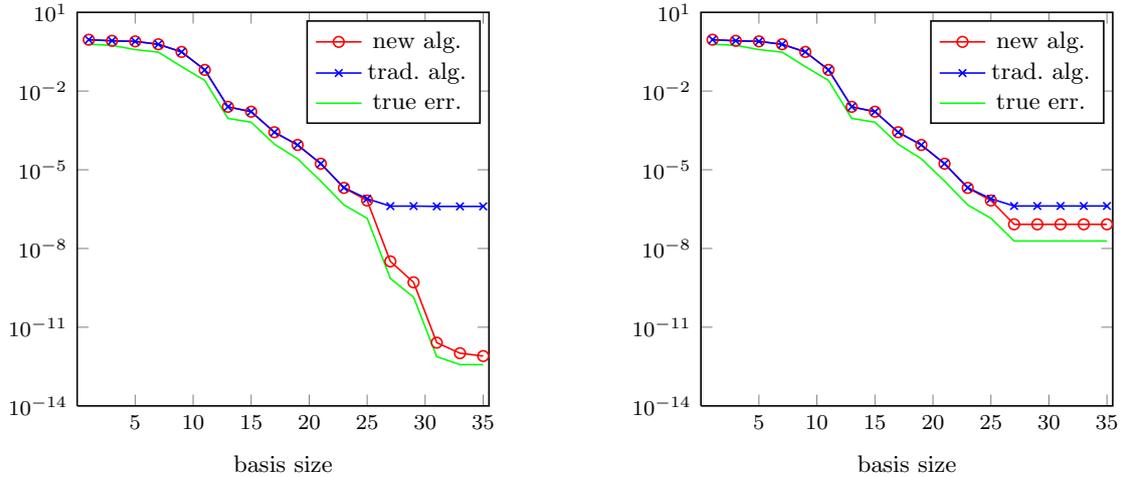
\begin{figure}
\begin{subfigure}[t]{0.45\textwidth}
\begin{center}
	\begin{tikzpicture}
		\begin{semilogyaxis}[xmin=0, xmax=35.5, ymin=1e-14, ymax=10, xtick={5, 10, 15, 20, 25, 30, 35},
			semithick, small, width=7.0cm, height=6.8cm, 
			minor tick style={draw=none},
                                     xlabel={basis size}, legend style={font=\small}]
			\addplot [mark=o, red] table[x=N, y=est0] {errors-5.dat};
			\addplot [mark=x, blue] table[x=N, y=est1] {errors-5.dat};
			\addplot [mark=none, green] table[x=N, y=err] {errors-5.dat};
                        \legend{new alg., trad.~alg., true err.}
		\end{semilogyaxis}
	\end{tikzpicture}
\end{center}
\caption{New algorithm used for basis generation}
\label{fig:errors}
\end{subfigure}
\hfill
\begin{subfigure}[t]{0.45\textwidth}
\begin{center}
	\begin{tikzpicture}
		\begin{semilogyaxis}[xmin=0, xmax=35.5, ymin=1e-14, ymax=10, xtick={5, 10, 15, 20, 25, 30, 35},
			semithick, small, width=7.0cm, height=6.8cm, 
			minor tick style={draw=none},
                                     xlabel={basis size}, legend style={font=\small}]
			\addplot [mark=o, red] table[x=N, y=est0] {errors-std-5.dat};
			\addplot [mark=x, blue] table[x=N, y=est1] {errors-std-5.dat};
			\addplot [mark=none, green] table[x=N, y=err] {errors-std-5.dat};
                        \legend{new alg., trad.~alg., true err.}
		\end{semilogyaxis}
	\end{tikzpicture}
\end{center}
\caption{Traditional algorithm used for basis generation}
\label{fig:traderrors}
\end{subfigure}
	\caption{Maximum relative reduction errors and estimated reduction
		errors ($H^1$-norm) for numerical example
		(\ref{eq:thermalblock}).}
\end{figure}

Equation (\ref{eq:thermalblock}) was discretized using linear finite elements on a regular grid with
$500\times 500 \times 2$ triangular entities (Fig.~\ref{fig:solution}). 
Then, a reduced space of dimension 35 was generated with the weak greedy algorithm
using our new algorithm for the evaluation of the error estimator.
An equidistant training set of $5^4$ parameters was used.
Finally, for each $n$-dimensional reduced subspace $\tilde{V}_n$ ($0\leq n \leq
35)$ produced by the greedy algorithm we computed the maximum reduction error and the maximum estimated reduction errors
using both the traditional and our improved algorithm on 20 randomly selected new parameters in
$\mathcal{P}$ (Fig.~\ref{fig:errors}).
Moreover, the maximum and minimum efficiencies (i.e.~the quotient error/estimate) of the estimator evaluated using both 
algorithms
were determined for the same random parameters (Table \ref{tab:efficiencies}). Our results clearly indicate the breakdown
of the traditional algorithm for more than 25 basis vectors at a relative error of about $10^{-7} \approx
\sqrt{\varepsilon}$ whereas our new algorithm remains efficient for all tested basis sizes.

To underline the need for accurate error estimation in order to obtain reduced spaces of high 
approximation quality, we repeated the same experiment using the traditional algorithm for error
estimation during basis generation (Fig. \ref{fig:traderrors}). While the maximum model reduction error
still improves from $10^{-7}$ to $10^{-8}$ after the breakdown of the error estimator,
the final reduced space approximates the solution manifold 4 orders of magnitude worse than
the space obtained with our improved algorithm.

\subsection*{Acknowledgements}
This work has been supported by the German Federal Ministry of Education and Research (BMBF) under contract number 05M13PMA
and by CST - Computer Simulation Technology AG.

\begin{table}
\caption{Maximum and minimum efficiencies ($H^1$-norm) of traditional and new error estimator for
	numerical example (\ref{eq:thermalblock}); efficiencies were calculated for 20 randomly chosen parameters. }
\label{tab:efficiencies}
\begin{center}
\setlength{\tabcolsep}{2.5mm}
\begin{tabular}{ llcccccc }
\toprule
\multicolumn{2}{r}{basis size} & 10&15&20&25&30&35 \\
\midrule
trad. & max & $4.9\cdot 10^{-1}$&$4.3\cdot 10^{-1}$&$4.6\cdot 10^{-1}$&$4.1\cdot 10^{-2}$&$1.9\cdot 10^{-5}$&$9.0\cdot 10^{-7}$ \\
      & min & $2.1\cdot 10^{-1}$&$2.3\cdot 10^{-1}$&$1.8\cdot 10^{-1}$&$1.0\cdot 10^{-3}$&$3.1\cdot 10^{-6}$&$4.8\cdot 10^{-7}$ \\
new   & max & $4.9\cdot 10^{-1}$&$4.3\cdot 10^{-1}$&$4.7\cdot 10^{-1}$&$3.9\cdot 10^{-1}$&$4.6\cdot 10^{-1}$&$4.6\cdot 10^{-1}$ \\
      & min & $2.1\cdot 10^{-1}$&$2.3\cdot 10^{-1}$&$2.2\cdot 10^{-1}$&$2.1\cdot 10^{-1}$&$2.3\cdot 10^{-1}$&$2.4\cdot 10^{-1}$ \\
\bottomrule
\end{tabular}
\end{center}
\end{table}

\FloatBarrier


\begin{thebibliography}{99}
\setlength{\parskip}{0pt}

\bibitem{RBBook} A. T. Patera, G. Rozza. {\it Reduced basis approximation and a posteriori error estimation for parametrized partial differential equations}, Version 1.0, Copyright MIT 2006, to appear in (tentative rubric) MIT Pappalardo Graduate Monographs in Mechanical Engineering.
\bibitem{Haasdonk} B. Haasdonk, M. Ohlberger. Reduced basis method for finite volume approximations of parametrized linear evolution equations. {\em M2AN (Math. Model. Numer. Anal.)}, Vol. {\bf 42(2)}, 277--302, 2008.
\bibitem{Binev} P. Binev, A. Cohen, W. Dahmen, R. DeVore, G. Petrova, P. Wojtaszczyk. Convergence Rates for Greedy Algorithms in Reduced Basis Methods. {\em SIAM J. Math. Anal.}, Vol. {\bf 43(3)}, 1457--1472, 2011.
\bibitem{YANO} M. Yano. A space-time Petrov-Galerkin certified reduced basis method: Application to the boussinesq equations.  Accepted in {\em SIAM Journal on Scientific Computing}, 2013.
\bibitem{BennerHess} P. Benner, M. Hess. The Reduced Basis Method for Time-Harmonic Maxwell's Equations. {\em Proceedings in Applied Mathematics and Mechanics}, Vol. {\bf 12}, 661--662, 2012.
\bibitem{Casenave12} F. Casenave. Accurate a posteriori error evaluation in the reduced basis method. {\em C. R. Math. Acad. Sci}, Vol. {\bf 350}, 539--542, 2012.
\bibitem{Casenave13} F. Casenave, A. Ern, T. Leli\`{e}vre. Accurate and online-efficient evaluation of the a posteriori error bound in the reduced basis method. Accepted in {\em M2AN (Math. Model. Numer. Anal.)}, 2013.
\end{thebibliography}
\end{document}